Zeilberger
 
\def\Re{{\rm Re}\,}
\def\halmos{\hbox{\vrule height0.15cm width0.01cm\vbox{\hrule
height
 0.01cm width0.2cm \vskip0.15cm \hrule height 0.01cm
width0.2cm}\vrule
 height0.15cm width 0.01cm}}
\font\eightrm=cmr8  
\font\eighttt=cmtt8
\magnification=\magstephalf

\parindent=0pt
\overfullrule=0in
\bf
\noindent
A HIGH-SCHOOL ALGEBRA\footnote{$^1$}
{\eightrm 
and high-school (purely formal) calculus.
}
, WALLET-SIZED PROOF, OF THE BIEBERBACH 
\hbox{ \hskip96pt \relax CONJECTURE [After L. Weinstein]}
\medskip
\noindent
\it
\hskip100pt \relax Shalosh B. Ekhad${}^2$ and Doron
Zeilberger\footnote{$^2$}
{\eightrm 
Department of Mathematics, Temple University,
Philadelphia, PA 19122, USA. 
\break
{\eighttt ekhad@euclid.math.temple.edu ,}
{\eighttt zeilberg@euclid.math.temple.edu .}
Supported in part by the NSF. We would like to thank Richard Askey
and Jeff Lagarias for comments, that improved readability.
}
\bigskip
\it
\qquad \qquad \qquad \qquad \qquad \qquad
\qquad\qquad Dedicated to Leonard Carlitz\footnote{$^3$}
{\eightrm
L. Carlitz was, for many years, editor of the Duke Journal, until
he was relieved from his duties by the proponents of so-called
``modern math'', who proceeded to reject anything that smacked,
even
faintly, of Carlitz-style mathematics. But those who have drowned
Carlitz will soon be drowned themselves, as post-modern, 
computer-assisted and computer-generated mathematics, that by its
very nature
is purely formal, will soon take over and make so-called
``modern'' mathematics a relic of the past.
}
, master of formal mathematics
\bigskip
\rm
Weinstein's[2]  brilliant short proof of de Branges'[1] theorem can
be made
yet much shorter(modulo routine calculations), completely
elementary
(modulo L\"owner theory),  self contained(no need for the esoteric
Legendre polynomials' addition theorem), and motivated(ditto),
as follows. Replace the text between p. 62, line 7 and p. 63, line
7, by Fact 1
below, and
the text between the last line of p.63 and p.64, line 7, by Fact 2
below.
\par
\smallskip
\noindent
{\bf FACT 1:} Let $f_t (z) = e^t z \exp ( \sum_{k=0}^{\infty} c_k
(t) z^k )$ 
where $c_k (t)$ are {\it formal} functions
of $t$. Let $z$ and $w$ be related by $z/(1-z)^2 = e^t w/(1-w)^2$.
The 
following formal identity holds. (For any formal Laurent series
$f(z)$,
$CT_z f(z)$ denotes the {\it Constant Term} of $f(z)$.)
\smallskip
$$
(1+w) { {d} \over {dt}} \{
\sum_{k=1}^{\infty} ( 4/k - k  c_k (t) \overline{c_k (t)} ) w^k \}
=
$$
$$
(1-w) \sum_{k=1}^{\infty} \Re CT_z
\left\{
{
{ {\partial f_t (z)} \over {\partial t} }
\over
{ z {\partial f_t (z)} \over {\partial z} } 
}
\cdot
(2(1+  \dots + k c_k (t) z^k ) - k c_k (t) z^k ) \cdot
(2(1+  \dots + k \overline{c_k (t}) z^{-k} ) - k \overline{c_k (t)}
z^{-k} ) 
\right\} w^k
$$
\medskip
\noindent
{\bf Proof:} Routine. (Obviously computer-implementable.) \halmos
\bigskip
\noindent
{\bf FACT 2:} The polynomials $A_{k,n} (c)$, defined in terms of
the
formal power series (Laurent in $w$) expansion
$(1- z(2c+ (1-c)(w+ 1/w))+z^2 )^{-1}=
\sum_{n=0}^{\infty}\sum_{k=0}^{n}
A_{k,n}(c) (w^k + w^{-k}) z^n$ are non-negative.
\medskip
\noindent
{\bf Proof:} This follows immediately from the stronger fact that
the 
polynomials$B_{k,n} (c)$, defined by the expansion
$(1- z(2c+ (1-c)(w+ 1/w))+z^2
)^{-1/2}=\sum_{n=0}^{\infty}\sum_{k=0}^{n}
B_{k,n}(c) (w^k + w^{-k}) z^n$ are perfect squares. To prove that
\smallskip
$$
B_{k,n} (c) := CT_{z,w} F(z,w,c,k,n) := CT_{z,w}
 [ {
{(1- z(2c+ (1-c)(w+ 1/w))+z^2 )^{-1/2}}
\over { z^{n} w^{k}}
} ]
$$
\smallskip
are indeed perfect squares,
the reader\footnote{$^4$}
{\eightrm
Human readers: Find a computer friend to help you with this. All
you have to do is express $G_1$, $G_2$, qua polynomials in $z,w$
(of degree $2$ in each), generically, with "indeterminate
coefficients",
then divide (WZ) by $F$, simplify
clear denominators, and equate all the coefficients
of the monomials $z^i w^j$ in the resulting identity to $0$,
getting
a linear system of equations, with 2[(2+1)$\cdot$(2+1)]+4
unknowns, that your computer friend can easily solve. This method
is
called the WZ method, and the fact that such a recurrence {\it
always}
exists follows from the WZ theory (Wilf and Zeilberger,
Invent. Math.108(1992), 575-633.), but at any particular instance,
like in this case, no explicit reference to WZ theory need be made.
}
can easily find polynomials in $(n,k,c)$,
$p_0,p_1,p_2,p_3$, and polynomials in $(n,k,c,z,w)$, $G_1$, and
$G_2$, 
both of degree $2$ in both $z$ and $w$, such that
\smallskip
$$
p_0 F(z,w,c,k,n) +  p_1 F(z,w,c,k,n+1) + p_2 F(z,w,c,k,n+2) +
p_3 F(z,w,c,k, n+3 ) =
$$
$$
z {{d} \over {dz}}  ( {{G_1 F} \over {z^3 w}}) +
w {{d} \over {dw}}   ( {{G_2 F} \over {z w^3}}) \quad .
\eqno(WZ)
$$
\smallskip
Applying $CT$ to both sides of (WZ), remembering the obvious fact
that for {\it any} formal Laurent series $f(z)$, $ CT (
z(d/dz)f(z))=0$,
we get that the $B_{k,n}$ satisfy the linear recurrence, in $n$:
\smallskip
$$
p_0 B_{k,n}(c) +  p_1 B_{k,n+1}(c) + p_2 B_{k,n+2}(c)+
p_3 B_{k,n+3} (c) =0 \qquad.
\eqno(Rec^{2}) 
$$
\smallskip
The recurrence $(Rec^{2})$ can be used to generate many
$B_{k,n}(c)$,
and it turns out, empirically for now, that they are all perfect
squares $B_{k,n} (c) = L_{k,n} (c)^2$, for some double-sequence of
polynomials $L_{k,n}(c)$. These empirically-generated polynomials
can be used to find a (conjectured) linear recurrence
\smallskip
$$
q_0 L_{k,n} (c) + q_1 L_{k,n+1} (c) + q_2 L_{k,n+2}(c)
=0 \quad ,
\eqno(Rec)
$$
\smallskip
where $q_0 , q_1 , q_2$ are polynomials of $(n,k,c)$. Let's now
{\it define} $L'_{k,n}(c)$ to be the solution of (Rec) under the
appropriate
initial conditions $L_{k,0}$, $L_{k,1}$, and define 
$B'_{k,n}(c) := L'_{k,n} (c)^2$. Using high-school linear algebra
(which is implemented in the {\it gfun} Maple package developed
by Salvy and Zimmerman) one can easily find a (third order)
recurrence satisfied by $B'_{k,n}(c)$, that turns out to be
identical with
($Rec^2$).  Matching the three initial values $n=0,1,2$ 
completes the proof. \halmos
\medskip
\noindent
{\bf References}
\smallskip
\noindent
1. L. de Branges, {\it A proof of the Bieberbach conjecture}, Acta
Math.
{\bf 154}(1985), 137-152.
\par
\noindent
2. L. Weinstein, {\it  The Bieberbach Conjecture}, Duke Math. J.
{\bf 59}(1991)
61-64.
\bye